\documentclass[10pt,leqno]{article}
\usepackage{amssymb,amsbsy,amsmath,amsfonts,amssymb,amscd, mathrsfs}

\usepackage[english]{babel}
\usepackage[T1]{fontenc}
\usepackage{indentfirst}

\makeatletter
\@addtoreset{equation}{section}
\makeatother

\newtheorem{statement}{}[section]
\newtheorem{theoreme}[statement]{Theorem}
\newtheorem{lemme}[statement]{Lemma}

\newcommand\C{\mathbb C}

\newcommand\R{\mathbb R}
\newcommand\T{\mathbb T}
\newcommand\D{\mathbb D}

\newcommand\e{{\rm e}}

\newcommand\eps{\varepsilon}
\newcommand\ind{{\rm 1\kern-.30em I}}
\newcommand\qed{\hfill $\square$}
\renewcommand \Re{{\mathfrak R}{\rm e}\,}

\let\phi=\varphi

\title{\bf Nevanlinna counting function and Carleson function of analytic maps}
\author{\it Pascal Lef\`evre, Daniel Li,\\ \it Herv\'e Queff\'elec, Luis Rodr{\'\i}guez-Piazza}

\date{\footnotesize \today}

\begin{document}

\maketitle

\noindent{\bf Abstract.} \emph{We show that the maximal Nevanlinna counting function and the Carleson 
function of analytic self-maps of the unit disk are equivalent, up to constants.} 
\medskip

\noindent{\bf Mathematics Subject Classification.} Primary: 30C80 -- Secondary: 47B33; 47B10
\medskip

\noindent{\bf Key-words.}  Analytic self-map of the unit disk -- Carleson function -- Carleson measure -- 
composition operator -- Nevanlinna counting function 


\section{Introduction}

Carleson measures and the Nevanlinna counting function are two classical concepts in Complex Analysis. 
Carleson measures emerged in 1958 when L. Carleson (\cite{Carleson1}, \cite{Carleson2}) showed his famous 
embedding theorem: For any positive finite measure $\mu$ on the closed unit disk $\overline \D$, 
the identity map from the Hardy space $H^2$ into $L^2 (\mu)$ is bounded if and only if this measure 
satisfies the following geometric condition: $\sup_{|\xi| = 1} \mu [ W (\xi, h) ] = O\,(h)$, where $W(\xi, h)$ is the 
Carleson window of size $h$ centered at $\xi$. This supremum is called the Carleson function $\rho_\mu$ 
of $\mu$. \par
If $\phi$ is an analytic self-map of $\D$ (such a function is sometimes called a Schur function), $\phi$ induces a 
composition operator $C_\phi \colon f \in H^2 \mapsto f\circ \phi \in H^2$, which may be seen as the identity 
from $H^2$ into $L^2 (m_\phi)$, where $\mu = m_\phi$ is the image of the Lebesgue measure on the unit circle 
by $\phi^\ast$, the boundary values function of $\phi$. We say that $\rho_\phi = \rho_{m_\phi}$ is the Carleson 
function of $\phi$.\par
\medskip

Nevanlinna counting function traces back earlier, in the thirties of the last century, in connection with the Jensen 
formula and the Nevanlinna theory of defect (\cite{Rubel-livre} or \cite{Nevanlinna-livre}). It is defined, for 
$w \in \phi (\D)$ and $w \neq \phi (0)$, by $N_\phi (w) = \sum_{\phi (z) = w} \log 1/|z|$ (see 
\eqref{definition Nevanlinna}).\par
In a slightly different context,  Littlewood used it implicitly (\cite{Littlewood}, see Theorem~4) when he showed 
that, for every analytic self-map $\phi$ of $\D$, we have $N_{\varphi}(z) = O\,(1 - | z |)$ as $| z | \to 1$. This 
turns out to imply (\cite{Shap}, \cite{Shap-livre}) that the composition operator 
$f \mapsto f \circ \phi = C_{\varphi}(f)$ is continuous on $H^2$ (which precisely  means, in present language, 
that $m_\phi$ is a Carleson measure).
\par\bigskip

Later, and till now, the regularity of  composition operators $C_\varphi$ on $H^2$ (their compactness, or 
membership in a Schatten class) in terms of their ``symbol'' $\varphi$ has been studied either from the point 
of view of Carleson measures or from the point of view of the Nevanlinna counting function, those two points of 
view being completely  separated. For example, the compactness of $C_\varphi \colon H^2\to H^2$ has been 
characterized in terms of the Carleson function of the symbol $\rho_\phi (h) = o\, (h)$, as $h \to 0$, by 
B. McCluer (\cite{McCluer} -- see also \cite{Power}). In another paper, it was characterized  in terms of the 
Nevanlinna counting function $N_\phi$ of the symbol: $N_\phi (w) = o\, (1 - |w|)$, as $|w| \to 1$, by 
J. Shapiro (\cite{Shap}). A similar situation exists for the characterization of the  membership of $C_\varphi$ 
in a prescribed  Schatten class (\cite{Luecking} and \cite{Luecking-Zhu}).\par
\medskip 

Though the definition the Carleson measure $m_\varphi$ and that of the Nevanlinna counting function 
$N_\varphi$ are of different nature, there should therefore exist a direct link between these two quantities.\par
\medskip

Some results in this direction had been given: B. R. Choe (\cite{Choe}) showed that 
$\limsup_{h \to 0} (\rho_\phi (h) / h)^{1/2} $ is equivalent, up to constants, to the distance of $C_\phi$ to 
the space of compact operators on $H^2$; since J. Shapiro proved (\cite{Shap}) that this distance is 
$\limsup_{|w| \to 1} (N_\phi (w) / \log |w|)^{1/2}$, one gets that 
\begin{displaymath}
\limsup_{|w| \to 1} N_\phi (w) / \log |w| \approx \limsup_{h \to 0} \rho_\phi (h) / h\,.
\end{displaymath}
Later,  J. S. Choa and H. O. Kim (\cite{Choa-Kim}) gave a somewhat direct proof of the equivalence of the two above 
conditions, without using the properties of the composition operator, but without giving explicitly a direct 
relation between the two functions $\rho_\phi$ and $N_\phi$. \par
\bigskip

The aim of this paper is to show the surprising fact that the Nevanlinna counting function and the Carleson function 
are actually equivalent, in the following sense: 
\par
\bigskip

\begin{theoreme}\label{theo principal}
There exists a universal constant $C > 1$, such that, for every analytic self-map $\phi \colon \D \to \D$, one has:
\begin{equation}
(1/C) \, \rho_\phi (h/C) \leq \sup_{|w| \geq 1 - h} N_\phi (w) \leq C\,\rho_\phi (C\,h),
\end{equation} 
for $0 < h < 1$ small enough.\par
More precisely, for every $\xi \in \partial \D$, one has:
\begin{equation}
(1/64) \, m_\phi [ W (\xi, h/64) ] \leq \sup_{w \in W (\xi, h) \cap \D} N_\phi (w) 
\leq 196\, \, m_\phi [ W (\xi, 24\,h)] \,,
\end{equation}
for $0 < h < (1 - |\phi (0)|) / 16 $.
\end{theoreme}
\bigskip

Actually the above explicit constants are not relevant and we did not try to have ``best'' constants. It can be 
shown that for every $\alpha > 1$, there is a constant $C_\alpha > 0$ such that
$m_\phi \big( S (\xi, h) \big) \leq C_\alpha \, \tilde\nu_\phi (\xi, \alpha h) $ 
and 
$\tilde\nu_\phi (\xi, h) \leq C_\alpha \, m_\phi \big( S (\xi, \alpha h) \big)$ for 
$0 < h < (1 - |\phi (0)|)/ \alpha$, where $S (\xi, h)$ is defined in \eqref{definition S} and 
$\tilde \nu (\xi, h) = \sup_{w \in S (\xi, h) \cap \D} N_\phi (w)$ (see \eqref{definition nu tilde}).

\section{Notation}

We shall denote by $\D = \{ z\in \C\,; \ |z| < 1\}$ the open unit disc of the complex plane and by 
$\T = \partial \D = \{ z\in \C\,;\ |z | = 1\}$ its boundary; $m$ will be the normalized Lebesgue 
measure $dt / 2\pi$ on $\T$, and $A$ the normalized Lebesgue measure $dxdy / \pi$ on 
$\overline \D$. 
For every analytic self-map $\phi$ of $\D$, $m_\phi$ will be the pull-back measure of $m$ by $\phi^\ast$, 
where $\phi^\ast$ is the boundary values function of $\phi$. 
\par
For every $\xi \in \T$ and $0 < h < 1$, the Carleson window $W (\xi, h)$ centered at $\xi$ and of size $h$ 
is the set
\begin{equation}
W (\xi, h) = \{ z \in \overline{\D}\,;\ |z| \geq 1 - h\quad \text{and} \quad |\arg ( z \bar{\xi})| \leq h\}.
\end{equation}
For convenience, we shall set $W (\xi, h) = \overline \D$ for $h \geq 1$.\par
For every analytic self-map $\phi$ of $\D$, one defines  the maximal function of $m_\phi$, for $0 < h < 1$, 
by:
\begin{equation}
\rho_\phi (h) = \sup_{\xi \in \T} m \big(\{\zeta \in \T\,;\ \phi^\ast (\zeta) \in W (\xi, h) \} \big) 
= \sup_{\xi \in \T} m_\phi \big( W (\xi, h) \big).
\end{equation}
We have $\rho_\phi (h) = 1$ for $h \geq 1$. We shall call this function $\rho_\phi$ the 
\emph{Carleson function} of $\phi$. For convenience, we shall often also use, instead of the Carleson 
window $W (\xi, h)$, the set 
\begin{equation}\label{definition S}
S (\xi, h) =\{ z\in \overline\D\,;\ |z - \xi | \leq h\}\,,
\end{equation}
which has an equivalent size.
\par
\smallskip

The Nevanlinna counting function $N_\phi$ is defined, for $w \in \phi(\D) \setminus \{\phi (0) \}$, by
\begin{equation}\label{definition Nevanlinna}
N_\phi (w) = \sum_{\phi (z) = w} \log \frac{1}{|z|}\,,
\end{equation}
each term $\log \frac{1}{|z|}$ being repeated according to the multiplicity of $z$, and $N_\phi (w) = 0$ 
for the other $w \in \D$. Its maximal function will be denoted by
\begin{equation}\label{definition nu phi}
\nu_\phi (t) = \sup_{|w| \geq 1 - t} N_\phi (w).
\end{equation}
%


\section{Majorizing the Nevanlinna counting function by the Carleson function}

The goal of this section is to prove:

\begin{theoreme}\label{th majoration Nevanlinna}
For every analytic self-map $\phi$ of $\D$, one has, for every $a \in \D$:
\begin{equation}\label{eq majoration Nevanlinna}
N_\phi (a) \leq 196\, m_\phi \big( W(\xi, 12 h) \big),
\end{equation}
for $0 < h < (1 - |\phi (0)|) / 4 $, where $\xi =\frac{a}{|a|}$ and $h = 1 - |a|$.\par
In particular, for $0 < h < (1 - |\phi (0)|) / 4 $:
\begin{equation}
\nu_\phi (h) = \sup_{|a| \geq 1 - h} N_\phi (a) \leq 196\,\rho_\phi (12 h).
\end{equation}
\end{theoreme}

Let us note that, since $W (\zeta, s) \subseteq W (\xi, 2t)$ whenever $0 < s \leq t$ and 
$\zeta \in W (\xi, t) \cap \partial \D$, we get from 
\eqref{eq majoration Nevanlinna} that
\begin{equation}
\sup_{w \in W (\xi, h) \cap \D} N_\phi (w) \leq 196\, m_\phi \big( W(\xi, 24 h) \big) \,.
\end{equation}

\bigskip

We shall first  prove the following lemma.
\begin{lemme}\label{lemme preparatoire}
Let $\phi$ be an analytic self map of $\D$. For every $z\in \D$, one has, if $w = \phi (z)$, $\xi = w / |w|$ and 
$h = 1 - |w| \leq 1/4$:
\begin{equation}
m_\phi \big( W (\xi,  12\, h) \big) \geq m_\phi \big( S (\xi, 6h) \big) 
\geq \frac{|w|}{8}\, (1 - |z|)\,.
\end{equation}
\end{lemme}

\noindent{\bf Proof.} 
We may assume, by making a rotation, that $w$ is real and positive: $3/4 \leq w < 1$.\par
Let:
\begin{equation}\label{def T}
T (u) =\frac{a u + 1}{u + a}\,\raise 1 pt \hbox{,} 
\end{equation}
where
\begin{displaymath}
\quad a = w - \frac{2}{w} < - 1\, ,
\end{displaymath}
so that $T \colon \D \to \D$ is analytic, and $T (w ) = w/2$.\par
If $P_z$ is the Poisson kernel at $z$, one has:
\begin{displaymath}
\frac{w}{2} = T[\phi (z)] = \int_\T (T\circ \phi)^\ast P_z\,dm =\int_\T \Re [(T\circ \phi)^\ast] P_z\,dm.
\end{displaymath}
Hence, if one sets:
\begin{displaymath}
E = \{ \Re (T \circ \phi^\ast) \geq w/4\} = \{ \Re [(T \circ \phi)^\ast] \geq w/4\},
\end{displaymath}
one has:
\begin{displaymath}
\frac{w }{2} \leq \int_E P_z \,dm + \frac{w}{4} \int_{E^c} P_z\,dm 
\leq \int_E P_z \,dm + \frac{w}{4} \int_{\D} P_z\,dm = \int_E P_z \,dm + \frac{w}{4} \,;
\end{displaymath}
therefore:
\begin{displaymath}
\int_E P_z\,dm \geq \frac{w}{4}\,\cdot
\end{displaymath}
Since 
\begin{displaymath}
\| P_z\|_\infty = \frac{1 + |z|}{1 - |z|} \leq \frac{2}{1 - |z|} \,\raise 1,5pt \hbox{,}
\end{displaymath}
we get:
\begin{equation} \label{minoration de m(E)}
m (E) \geq \frac{w}{8}\, (1 - |z|) \,.
\end{equation}
On the other hand, \eqref{def T} writes
\begin{equation}\label{inverse T}
u = T^{-1} (U) = \frac{aU - 1}{a - U}\,;
\end{equation}
hence:
\begin{displaymath}
|1 - u | = |a + 1| \frac{|1 - U|}{|a - U|} \leq  \frac{2\, |a + 1|}{|a - U|}\,\cdot
\end{displaymath}
But $a < - 1$ is negative, so $\Re U \geq w/4$ implies that 
\begin{displaymath}
|a - U| \geq \Re (U - a) \geq \frac{w}{4} - a = \frac{2}{w} - \frac{3}{4}\,w \geq \frac{5}{4} \,\cdot
\end{displaymath}
Moreover, for $w \geq 3/4$:
\begin{displaymath}
|a + 1| = (1 - w)\,\Big(\frac{2}{w}  + 1\Big) \leq  \frac{11}{3}\, (1 - w)\,.
\end{displaymath}
We get hence $| 1 - u| \leq 6\, h $ when \eqref{inverse T} holds and $\Re U \geq w/4$.\par
It follows that:
\begin{equation}
\phi^\ast (E) \subseteq T^{-1} (\{\Re U \geq w/4\}) \subseteq S (1, 6 h),
\end{equation}
giving $m_\phi \big( W (1, 12 h) \big) \geq m_\phi \big( S (1, 6 h)\big) \geq m(E)$.\par\smallskip

Combining this with \eqref{minoration de m(E)}, that finishes the proof. \qed
\bigskip

\noindent{\bf Remark.} Theorem~\ref{th majoration Nevanlinna} follows immediately when $\phi$ is 
univalent since then, for $|w| \geq 3/4$ and $\phi (z) = w$:
\begin{displaymath}
N_\phi (w) = \log \frac{1}{|z|} \approx (1 - |z|) 
\lesssim m_\phi \big( W (1, 12 h) \big) \,.
\end{displaymath}
\par\bigskip

When proving the equivalence between the conditions $\rho_\phi (h) = o\, (h)$, as $h\to 0$, and 
$N_\phi (w) = o\, (1 - |w|)$, as $|w| \to 1$, J. S. Choa and H. O. Kim proved 
(see \cite{Choa-Kim}, page 112) the following inequality, for every analytic self-map 
$\phi \colon \D \to \D$ and every $w \in \D$, close enough to $1$:
\begin{equation}\label{ineq Choa-Kim}
N_\phi (w) \leq \frac{(1 -|w|^2)^2}{8 |w|^2} \int_{\partial \D} \frac{1}{|1 -\bar{w} \phi (z)|^2}\,dm (z)\,.
\end{equation}

This result follows from an Hilbertian method, \emph{viz.} Littlewood-Paley's identity: 
\begin{equation}\label{Littlewood}
\| f \circ \phi \|_2^2 = | f \circ \phi (0)|^2  + 2 \int_\D |f ' (w) |^2 N_\phi (w)\,dA (w)
\end{equation}
for every $f \in H^2$. 
With \eqref{ineq Choa-Kim}, one cannot go beyond the order $2$; for instance, we can deduce from 
\eqref{ineq Choa-Kim} (see the proof of Theorem~\ref{th majoration Nevanlinna} below), that, for 
$0 < h \leq 1/2$:
\begin{equation}\label{ineq Choa-Kim revisitee}
\sup_{|w| = 1 - h} N_\phi (w) \lesssim h^2 \int_0^{1/h^2} \rho_\phi \Big(\frac{1}{\sqrt t}\Big)\,dt 
\lesssim h^2 + h^2 \int_h^1 \frac{\rho_\phi (u)}{u^3}\,du.
\end{equation}
This is of course interesting only when the second term in the last sum is at most of order $h^2$, so, when the 
integral is bounded. Nevertheless, this result suffices to show that Shapiro's criterion of compactness for 
$C_\phi \colon H^2 \to H^2$ is implied by McCluer's one. Moreover, when the pull-back measure 
$m_\phi$ is an $\alpha$-Carleson measure ({\it i.e.} $\rho_\phi (h) \leq C\, h^{\alpha}$ for some constant 
$C > 0$), with $1 \leq \alpha \leq 2$, we get 
\begin{displaymath}
N_\phi (w) \lesssim h^2 + h^2 \int_h^1 \frac{u^\alpha}{u^3}\,du \lesssim h^2 +h^2 h^{\alpha -2} 
\lesssim h^\alpha.
\end{displaymath}
Recall (\cite{JFA}, Corollary~3.2) that, when $m_\phi$ is an $\alpha$-Carleson measure, the composition 
operator $C_\phi$ is in the Schatten class $S_p$ on the Hardy space $H^2$, for every $p > 2/(\alpha  - 1)$, 
and that $m_\phi$ is $\alpha$-Carleson for every $\alpha \geq 1$ when $C_\phi \colon H^\Psi \to H^\Psi$ 
is compact, if $\Psi$ is an Orlicz function satisfying the growth condition $\Delta^2$ (\cite{JMAA}, 
Theorem~5.2).
\par
But \eqref{ineq Choa-Kim revisitee} does not suffice for the compactness of 
$C_\phi \colon H^\Psi \to H^\Psi$ on general Hardy-Orlicz spaces (see \cite{CompOrli} or \cite{LLQR-CRAS}).
\par\smallskip

In order to prove Theorem~\ref{th majoration Nevanlinna}, we shall replace the Littlewood-Paley identity, by a 
more general formula, deduced from Stanton's formula (see \cite{Stanton}, Theorem~2).

\begin{theoreme}[Stanton's formula] 
For every analytic self-map $\phi \colon \D \to \D$ and every subharmonic function $G \colon \D \to \R$, 
one has:
\begin{equation}
\lim_{r \uparrow 1} \int_{\partial \D} G [ \phi (r \xi) ]\, dm (\xi) = 
G [ \phi (0)] + \frac{1}{2} \int_\D \Delta\,G(w) N_\phi (w)\,dA (w),
\end{equation}
where $\Delta$ is the distributional Laplacian.
\end{theoreme}

\noindent{\bf Proof of Theorem~\ref{th majoration Nevanlinna}.} 
If $a \notin \phi (\D)$, one has $N_\phi (a) = 0$, and the result is trivial. We shall hence assume that 
$a \in \phi (\D)$.\par\smallskip

Let $\Phi \colon [0, \infty) \to [0, \infty)$ be an Orlicz function, that is a non-decreasing convex function 
such that $\Phi (0) = 0$ and $\Phi (\infty) = \infty$, and we assume that $\Phi '$ is also an Orlicz function. 
In other words, $\Phi ''$ is an arbitrary non-negative and non-decreasing function and 
$\Phi ' (x) = \int_0^x \Phi '' (t)\,dt$ and $\Phi (x) = \int_0^x \Phi ' (t)\,dt$.
\par\smallskip

Let now $f \colon \D \to \C$ be an analytic function. We have, outside the zeroes of $f$, in writing 
$\Delta \Phi (| f |) = 4\partial {\bar \partial}  \Phi (\sqrt{| f |^2})$:
\begin{equation}\label{Laplacien}
\Delta \Phi (| f |) = \bigg[ \Phi '' (| f |) + \frac{ \Phi ' (| f |)}{| f |} \bigg]\,|f ' |^2.
\end{equation}
\par

We shall only use here that:
\begin{equation}\label{Laplacien bis}
\Delta \Phi (| f |) \geq  \Phi '' (| f |)\,|f '|^2
\end{equation}
(this is a not too crude estimate, since, $\Phi '$ being an Orlicz function, $\Phi ''$ is non-negative and 
non-decreasing, and hence $\Phi' (x) = \int_0^x \Phi '' (t)\,dt  \leq x \Phi'' (x)$ and 
$\Phi ' (x) = \int_0^x \Phi '' (t)\,dt \geq \int_{x/2}^x \Phi '' (t)\,dt \geq (x/2)\, \Phi '' ( x/2)$).
\par\smallskip

Set now, for $a \in \D$:
\begin{equation}
\qquad f_a (z) = \frac{1 - | a |}{1 - {\bar a} z} \,\raise 1pt \hbox{,}\quad z\in \overline\D.
\end{equation}
Since $\Phi (| f_a |)$ is subharmonic ($\Phi$ being convex and non-decreasing) and bounded, we can 
use Stanton's formula as:
\begin{equation}
\int_{\partial \D} \Phi (| f_a \circ \phi |)\,dm 
\geq \frac{1}{2} \int_\D \Phi'' (| f_a |)\,| f'_a |^2\,N_\phi\,dA.
\end{equation}
Let $h = 1 - | a |$. For $|z - a| < h$, one has
\begin{displaymath}
| 1 - {\bar a} z| = | (1 - |a|^2) + {\bar a} ( a - z)| \leq (1 - |a|^2) + | a - z| \leq 2h + h = 3h;
\end{displaymath}
Hence $| f_a (z) | \geq \frac{h}{3h} = \frac{1}{3}$ for $|z - a|< h$. It follows, since $\Phi ''$ is non-decreasing:
\begin{equation}
\int_{\partial \D} \Phi (| f_a \circ \phi |)\,dm 
\geq \frac{1}{2}\,\Phi '' \Big( \frac{1}{3}\Big) \int_{D(a, h)} | f'_a |^2\,N_\phi\,dA.
\end{equation}
Now, if $\phi_a (z) =\frac{a - z}{1 - {\bar a} z}$, one has 
$| f '_a (z) | = \frac{| a |}{1 + | a |}\, |\phi'_a (z)| \geq \frac{3}{7} \,|\phi'_a (z)|$ 
(we may, and do, assume that $1 - |a| = h \leq 1/4$); hence:
\begin{align*}
\int_{\partial \D} \Phi (| f_a \circ \phi |)\,dm 
& \geq \frac{1}{2}\,\Phi '' \Big( \frac{1}{3} \Big) \,\frac{9}{49} \int_{D(a, h)} | \phi'_a |^2\,N_\phi\,dA \\
& = \frac{9}{98} \, \Phi '' \Big( \frac{1}{3} \Big) \int_{\phi_a ( D(a, h) )} N_{\phi_a \circ \phi}\,dA
\end{align*}
(because $N_{\phi_a \circ \phi} \big( \phi_a (w) \big) = N_\phi (w)$ and $\phi_a^{-1} = \phi_a$).\par 
But $\phi_a \big( D(a, h) \big) \supseteq D (0, 1/3)$: indeed, if $|w| < 1/3$, then $w = \phi_a (z)$, with 
\begin{displaymath}
|a - z| = \bigg| \frac{(1 -|a|^2) w}{1 -{\bar a} w} \bigg| \leq (1 -|a|^2) \frac{|w|}{1 -|w|} 
< 2h \frac{1/3}{1 - 1/3} = h.
\end{displaymath}

We are going now to use the sub-averaging property of the Nevanlinna function 
(\cite{Shap-livre}, page 190, \cite{Shap}, \S~4.6, or \cite{Zhu}, Proposition~10.2.4): for every analytic self-map 
$\psi \colon \D \to \D$, one has
\begin{displaymath}
N_\psi (w_0) \leq \frac{1}{A (\Delta)} \int_\Delta N_\psi (w)\,dA (w) \,,
\end{displaymath}
for every disk $\Delta$ of center $w_0$ which does not contain $\psi (0)$.\par

This will be possible thanks to the following:

\begin{lemme}
For $1 - |a| < (1 - |\phi (0)|) / 4$, one has $|(\phi_a \circ \phi) (0) | > 1/3$.
\end{lemme}

\noindent{\bf Proof.} One has 
$| 1 - {\bar a}\, \phi (0) | \leq (1 - |a|^2) + |\bar a|\, |a - \phi (0)| \leq (1 - |a|^2) + |a - \phi (0)|$; 
hence:
\begin{align*}
| \phi_a \big( \phi (0) \big) | 
\geq \frac{|a - \phi (0)|}{(1 - |a|^2) + |a - \phi (0)|} 
& \geq  1 - \frac{1 - |a|^2}{(1 - |a|^2) + |a - \phi (0)|} \\
& \geq 1 - \frac{1 - |a|^2}{|a - \phi (0)|} 
\geq 1 - 2 \, \frac{1 - |a|}{|a - \phi (0)|} \,\cdot
\end{align*} 
But when $1 -|a| < (1 - |\phi (0)|)/4$, one has:
\begin{displaymath}
|a - \phi (0) | \geq |a| - |\phi (0)| = (1 - |\phi (0)|) - (1 - |a|) > 3 (1 - |a|)\,,
\end{displaymath}
and the result follows. \qed
\medskip

Hence:
\begin{displaymath}
\int_{D (0, 1/3)} N_{\phi_a \circ \phi}\,dA \geq \frac{1}{9}\, N_{\phi_a \circ \phi} (0) 
= \frac{1}{9} \, N_\phi (a),
\end{displaymath}
and
\begin{equation}\label{mino integrale Phi}
\int_{\partial \D} \Phi (| f_a \circ \phi |)\,dm \geq 
\frac{1}{98}\, \Phi '' \Big( \frac{1}{3} \Big)\,N_\phi (a).
\end{equation} 
\smallskip

We now have to estimate from above $\int_{\partial \D} \Phi (| f_a \circ \phi |)\,dm$. For 
that, we shall use the following easy lemma.

\begin{lemme}\label{lemme nombres complexes}
For every $\xi \in \partial \D$ and every $h \in (0,1/2]$, one has:
\begin{equation}\label{ineg nombres complexes}
\hskip 1,5 cm |1 - \bar{a} z|^2 \geq \frac{1}{4} (h^2 + |z - \xi|^2)\,, \qquad \forall z \in \overline\D,
\end{equation}
where $a = (1 - h) \xi$.
\end{lemme}

\noindent{\bf Proof.} The result is rotation-invariant; so we may assume that $\xi =1$ (and hence $a >0$). 
Write $z = 1 - r \e^{i \theta}$. Since $|z| \leq 1$ if and only if $r \leq 2 \cos \theta$, one has 
$\cos \theta \geq 0$ and hence $|\theta | \leq \pi/2$. Then:
\begin{align*}
\hskip 30pt | 1 - \bar{a} z |^2 
& = | 1 - a (1 - r \e^{i \theta})|^2 = | 1 - a + a r \e^{i\theta}|^2 \\
& = (1 - a)^2 +a^2 r^2 +2ar (1 - a) \cos \theta \\
& \geq (1 - a)^2 +a^2 r^2 \geq \frac{1}{4} (h^2 +r^2) = \frac{1}{4} (h^2 +|z -1|^2). 
\hskip 37,5pt \hbox{$\square$}
\end{align*}

Then:
\begin{align*}
\int_{\partial \D} \Phi (| f_a \circ \phi |)\,dm 
& = \int_{\overline{\D}} \Phi \bigg( \frac{1 -|a|}{| 1 - {\bar a} z|}\bigg)\, dm_\phi (z)\\
& \leq \int_{\overline{\D}} \Phi \bigg( \frac{2h}{(h^2 +|z - \xi|^2)^{1/2}}\bigg)\, dm_\phi (z),   
\quad \text{by \eqref{ineg nombres complexes} } \\
& = \int_0^{+\infty} 
m_\phi \Big( \Phi\Big( \frac{2h}{(h^2 +|z - \xi|^2)^{1/2}} \Big) \geq t \Big)\,dt \\
& =  \int_0^{+\infty} 
m_\phi \big( (h^2 +|z - \xi|^2)^{1/2} \leq 2h / \Phi^{-1} (t) \big)\, dt \\
& = \int_0^{\Phi (2)} 
m_\phi \big( (h^2 +|z - \xi|^2)^{1/2} \leq 2h / \Phi^{-1} (t) \big) \,dt \,,
\end{align*}
since $h \leq (h^2 +|z - \xi|^2)^{1/2} \leq 2h / \Phi^{-1} (t)$ implies $t \leq \Phi (2)$. We get:
\begin{displaymath}
\int_{\partial \D} \Phi (| f_a \circ \phi |)\,dm 
\leq \int_0^{\Phi (2)} m_\phi \big(|z - \xi| \leq 2h / \Phi^{-1} (t) \big) \,dt\,.
\end{displaymath}
We obtain from \eqref{mino integrale Phi}, by setting $u = 2h / \Phi^{-1} (t)$:
\begin{displaymath}
N_\phi (a) \leq \frac{98}{\Phi '' (1/ 3)} 
\int_{h}^\infty m_\phi \big( S (\xi, u) \big)\, \frac{2h}{u^2} \,
\Phi '\bigg( \frac{2h}{u} \bigg) \,du \,\cdot
\end{displaymath}
Since $\Phi ' (x) \leq x \Phi '' (x)$, we get:
\begin{equation}\label{majo Nevanlinna avant-derniere}
N_\phi (a) \leq \frac{98}{\Phi '' (1/ 3)} 
\int_{h}^\infty m_\phi \big( S (\xi, u) \big)\,\frac{4 h^2}{u^3}\, 
\Phi '' \bigg( \frac{2h}{u} \bigg)\,du.
\end{equation}
\par \smallskip

We are going now to choose suitably the Orlicz function $\Phi$. It suffices to define $\Phi '' $, for 
$a \in \D$ given (with $\xi = a / | a |$ and $h = 1 - |a| \leq 1/4$). By 
Lemma~\ref{lemme preparatoire}, since $a \in \phi (\D)$, there is a constant $c_0 > 0$, such that 
$m_\phi \big( S (\xi, c_0 h) \big) > 0$; we can hence set (note that $m_\phi \big( S (\xi, u) \big) \leq 1$): 
\begin{equation}
\Phi '' (v) = \left\{
\begin{array}{cl}
1 & \text{if } 0 \leq v \leq h \,, \\
\displaystyle \frac{1}{m_\phi \big( S (\xi, 2h /v ) \big)} 
& \displaystyle \text{if } h \leq v \leq 2/c_0 \,, \\
\displaystyle \frac{1}{m_\phi \big( S (\xi, c_0 h) \big)} 
& \displaystyle \text{if } v \geq 2/c_0 \,.
\end{array}
\right.
\end{equation}
It is a non-negative non-decreasing function, so the assumptions made on $\Phi$ at the beginning are satisfied. 
One has, since $m_\phi \big( S (\xi, u) \big) \,\Phi '' ( 2h /u ) \leq 1$:
\begin{displaymath}
\int_h^\infty m_\phi \big( S (\xi, u) \big) \,\frac{4 h^2}{u^3}\, 
\Phi '' \bigg( \frac{2 h}{u} \bigg)\,du 
\leq \int_h^\infty \frac{4 h^2}{u^3} \,du = 2.
\end{displaymath}
Since $c_0 \leq 6$, one has $h \leq 1 / 3 \leq 2/c_0 $ and hence 
$\Phi '' (1 / 3) = 1/ m_\phi \big( S (\xi, 6 h) \big)$; therefore 
\eqref {majo Nevanlinna avant-derniere} gives, for $h \leq (1 - |\phi (0) |) /4$:
\begin{equation}
N_\phi (a) \leq 196\, m_\phi \big( S (\xi, 6 h) \big),
\end{equation}
finishing the proof since $S (\xi, 6h) \subseteq W (\xi, 12 h)$.
\qed
\bigskip

\section{Domination of the Carleson function by the Nevanlinna function}

We cannot expect to estimate individually from above the $m_\phi$-measure of Carleson windows centered 
at $\xi = w/ |w|$ by $N_\phi (w)$, as in Theorem~\ref{th majoration Nevanlinna}. In fact, consider 
a conformal mapping $\phi$ from $\D$ onto $\D \setminus [0, 1[$. One has $N_\phi (t) = 0$ for every 
$t\in [0, 1[$, though $m_\phi \big( W (1, h) \big) > 0$ for every $h > 0$ (because 
$W (1, h) \supset W (\e^{ih/2}, h/2)$ and $m_\phi \big( W (\e^{ih/2}, h/2) \big) > 0$ by 
Lemma~\ref{lemme preparatoire}).\par
\smallskip

Let us give another example. Let $\phi (z) = (1 + z)/2$. Then:\par
\smallskip
a) One has $\phi (\e^{i\theta}) = (\cos \theta /2)\, \e^{i\theta/2}$ (with $ |\theta| \leq \pi$). Hence 
$\phi (\e^{i\theta})  \in W (\e^{i \theta_0}, h)$ if and only if $\cos (\theta/2) \geq 1 - h$ and 
$|(\theta/2) -\theta_0| \leq h$, \emph{i.e.} $2 (\theta_0 - h) \leq \theta \leq 2 (\theta_0 + h)$.\par 
Now, $1 - \cos (\theta/2) \leq \theta^2/8$, so the modulus condition is satisfied when $\theta^2 \leq 8h$; 
in particular when $|\theta| \leq 2 {\sqrt h}$.\par
For $\theta_0 = {\sqrt h}$, $m_\phi \big( W (\e^{i\theta_0}, h) \big)$ 
is bigger than the length of the interval 
\begin{displaymath}
[- 2 {\sqrt h}, 2 {\sqrt h} ] \cap [2 ({\sqrt h} - h), 2 ({\sqrt h} + h)] 
= [2{\sqrt h} - 2h , 2{\sqrt h}]\,,
\end{displaymath}
that is $2 h$. Therefore $m_\phi \big( W (\e^{i\theta_0}, h) \big) \geq 2 h$.
\par\smallskip

b) Let now $w = \phi (z)$. Write $w = \frac{1}{2} + r \,\e^{i\zeta}$ with $0 \leq r < 1/2$. Then, writing 
$r =  \frac{1}{2} - s$, one has $|z| = |2w - 1| = 2r$ and 
\begin{displaymath}
N_\phi (w) = \log \frac{1}{|z|} = \log \frac{1}{2r} = \log \frac{1}{1 -2s} \approx s.
\end{displaymath}
Now, $|w|^2 = \frac{1}{4} + r^2 + r \cos \zeta$ and 
\begin{displaymath}
h \approx 1 -|w|^2  = \frac{1}{2} ( 1 - \cos \zeta) + s (1 +\cos \zeta) - s^2 
\approx \frac{\zeta^2}{4} + 2s . 
\end{displaymath}
Writing $\zeta = s^{1/2\alpha}$, one gets:\par
$(i)\ $ for ``small'' $\zeta$ (\emph{i.e.} $0 < \alpha \leq 1$): $h  \approx s$, and so 
$N_\phi (w) \approx h$; \par 
$(ii)$ for ``large''  (\emph{i.e.} $\alpha \geq 1$): $h \approx s^{1/\alpha}$, and so  
$N_\phi (w) \approx h^\alpha$.
\medskip

On the other hand, $ w = \e^{i \zeta/2} [ (1 - s) \cos (\zeta/2) - i s \sin (\zeta/2) ]$; hence, when $s$ goes to 
$0$, one has 
\begin{displaymath}
\theta_w := \arg w = \frac{\zeta}{2} + \arctan \bigg[ \frac{s \sin (\zeta/2)}{(1 - s) \cos (\zeta/2)} \bigg] 
\sim \frac{\zeta}{2} \approx \zeta\,.
\end{displaymath}
For $\alpha \geq 1$, one has $h \approx s^{1/\alpha} = \zeta^2$, \emph{i.e.} $\zeta \approx \sqrt h$. Then, 
choosing $\alpha > 1$ such that $\zeta = \theta_0$, one has 
$m_\phi \big( W (w /|w|, h) \big) \approx h$, though $N_\phi (w) \approx h^\alpha \ll h$.\par\smallskip

One cannot hence dominate $m_\phi \big( W (w /|w|, h) \big)$ by $N_\phi (w)$.
\par\smallskip

We can remark that, nevertheless, in either case, one has $\rho_\phi (h) \approx h$ and 
$\nu_\phi (h) \approx h$.

\par\bigskip

We shall prove:

\begin{theoreme}\label{majoration mesure Carleson}
For every analytic self-map $\phi \colon \D \to \D$, one has, for every $\xi \in \partial\D$:
\begin{equation}\label{majoration de rho}
m_\phi \big( W (\xi, h) \big) \leq 64\, \sup_{w \in W (\xi, 64 h) \cap \D} N_\phi (w)\,,
\end{equation}
for $0 < h < (1 - |\phi (0)|)/ 16 $.
\end{theoreme}

\noindent{\bf Proof.} We shall set:
\begin{equation}
\nu_\phi (\xi, h) = \sup_{w \in W (\xi, h) \cap \D} N_\phi (w) \,.
\end{equation}

Note that
\begin{displaymath}
\nu_\phi (h) = \sup_{|\xi | = 1} \nu_\phi (\xi, h) \,,
\end{displaymath}
where $\nu_\phi$ is defined in \eqref{definition nu phi}

If  for some $h_0 > 0$, one has $\nu_\phi (\xi, h_0) = 0$, then  
$\phi (\D) \subseteq \D \setminus W (\xi, h_0)$, and hence $m_\phi \big( W (\xi, h)\big) = 0$ for 
$0 < h < h_0$. Therefore we shall assume that $\nu_\phi (\xi, h) > 0$. 
We may, and do, also assume that $h \leq 1/4$. By replacing $\phi$ by $\e^{i\theta}\phi$, it suffices to 
estimate $m_\phi \big( S(1, h)\big)$ (recall that $S (1, t) =\{ z\in \overline \D\,;\ |1 - z | \leq t \}$).\par
\smallskip

We shall use the same functions $f_a$ as in the proof of Theorem~\ref{th majoration Nevanlinna}, but, for 
convenience,  with a different notation. We set, for $0 < r < 1$:
\begin{equation}
u (z) = \frac{1 - r}{1 - rz} \,\cdot
\end{equation}
\par

Let us take an Orlicz function $\Phi$ as in the beginning of the proof of 
Theorem~\ref{th majoration Nevanlinna}, which will be precised later. We shall take this function in such 
a way that $\Phi \big( | u ( \phi (0) )| \big) = 0$.\par 
Since $\Phi ' (x) \leq x\, \Phi '' (x)$, \eqref{Laplacien} becomes:
\begin{equation}
\Delta \Phi (|u|) \leq 2 \Phi '' (|u|)\,|u'|^2, 
\end{equation}
and Stanton's formula writes, since $\Phi \big( | u ( \phi (0) )| \big) = 0$:
\begin{equation}\label{Stanton bis}
\int_{\partial\D} \Phi ( |u \circ \phi|) \,dm \leq 
\int_{\D} \Phi '' \big(|u (w)| \big)\,|u' (w)|^2\, N_\phi (w)\,dA (w).
\end{equation}
\par\smallskip

In all the sequel, we shall fix $h$, $0 < h \leq 1/4$, and take $r =  1 - h$.\par
\medskip

For $|z| \leq 1$ and $|1 - z| \leq h$, one has $|1 - r z|= |(1 - z) +hz | \leq |1 -z| + h \leq 2h$, so:
\begin{displaymath}
| u (z)| \geq \frac{(1- r)}{2h} = \frac{1}{2}\,\cdot
\end{displaymath}
Hence:
\begin{align*}
m_\phi\big( S(1,h) \big) 
& \leq \frac{1}{\Phi(1/2)} \int_{S(1, h)} \Phi \big( |u (z)| \big) \,dm_\phi (z) \\
& \leq  \frac{1}{\Phi(1/2)} \int_{\D} \Phi \big( |u (z)| \big) \,dm_\phi (z) \\
& =  \frac{1}{\Phi(1/2)} \int_\T  \Phi \big( |(u\circ \phi) (z)| \big) \,dm(z) \,,
\end{align*}
and so, by \eqref{Stanton bis}:
\begin{equation}\label{majoration avec Stanton bis}
m_\phi\big( S(1,h) \big)  \leq 
\frac{1}{\Phi(1/2)} \int_{\D} \Phi '' \big( |u (z)| \big) \, |u '(z)|^2 \, N_\phi (z)\,dA (z).
\end{equation}
\par

We are going to estimate this integral by separating two cases:  $| 1 - z| \leq h$ and $|1 - z| > h$.\par
For convenience, we shall set:
\begin{equation}\label{definition nu tilde}
\tilde \nu (t) = \sup_{w \in S (1, t) \cap \D } N_\phi (w) \,.
\end{equation}
\par\smallskip

1) Remark first that  
\begin{displaymath}
u ' (z) = \frac{r h}{( 1 - r z)^2}\,\raise 1,5pt \hbox{,}
\end{displaymath}
and so:
\begin{displaymath}
| u '(z)| \leq \frac{h}{(1 - r)^2} = \frac{1}{h}\,\cdot
\end{displaymath}
Since $| u (z) |\leq 1$, we get hence:
\begin{displaymath}
\int_{|1 - z| \leq h} \Phi'' \big(|u (z)| \big)\,|u '(z)|^2 N_\phi (z)\,dA(z) 
\leq \int_{S (1,h)} \Phi '' (1) \,\frac{1}{h^2}\, \tilde\nu (h)\,dA(z) \,, 
\end{displaymath}
giving, since $A \big( S (1, h) \big) \leq h^2$:
\begin{equation}\label{morceau 1}
\int_{|1 - z| \leq h} \Phi '' \big( |u (z)| \big)\,|u '(z)|^2 N_\phi (z)\,dA(z) 
\leq \Phi '' (1)\, \tilde\nu (h)\,.
\end{equation}
\par

2) For $0 < h \leq 1/4$, one has:
\begin{displaymath}
|u (z)| \leq \frac{2 h}{|1 - z|} \qquad \text{and} \quad | u ' (z) | \leq \frac{2 h}{\  | 1 - z|^2}\,;
\end{displaymath}
indeed, we have (this is obvious, by drawing a picture):
\begin{displaymath}
| 1 - r z| = r\,\Big|\frac{1}{r} - z\Big| \geq r\,|1 - z| \,,
\end{displaymath}
and hence $| 1 - r z| \geq \frac{3}{4}\, |1 - z|$, since $r = 1 - h \geq 3/4$.
We obtain:
\begin{displaymath}
\begin{array}{ll}
& \displaystyle \int_{|1 - z| > h} \Phi '' \big( |u (z)| \big)\, |u '(z)|^2 N_\phi (z)\,dA(z) \\
& \displaystyle \hskip 100 pt \leq 4 \int_{|1 - z| > h} \Phi '' \bigg( \frac{2h}{|1 - z|}\bigg)\, 
\frac{h^2}{|1 - z|^4} N_\phi (z)\,dA(z).
\end{array}
\end{displaymath}

Then, using polar coordinates centered at $1$ (note that we only have to integrate over an arc of length less 
than $\pi$), and the obvious inequality $N_\phi (z) \leq \tilde \nu (|1 - z|)$, we get:
\begin{align}\label{morceau 2}
& \int_{|1 - z| > h} \Phi '' \big( |u (z)| \big)\, |u '(z)|^2 N_\phi (z)\,dA(z) \\
& \displaystyle \hskip 120 pt 
\leq 4 \int_h^2 \Phi '' \bigg(\frac{2h}{t}\bigg)\, \frac{h^2}{t^3}\,\tilde\nu (t) \,dt \,. \notag
\end{align}
\par\smallskip

We now choose the Orlicz function as follows (with $a = \phi (0)$):
\begin{equation}
\Phi '' (v) = \left\{
\begin{array}{cl}
\displaystyle 0 & \text{if } 0 \leq v \leq h / ( 1 - |a|) \,, \\
\displaystyle \frac{1}{\tilde\nu (2 h/ v)} & \text{if } h / ( 1 - |a|) < v < 2\,,\\
\smallskip
\displaystyle \frac{1}{\tilde \nu (h)} & \text{if } v \geq 2 \,.
\end{array}
\right.
\end{equation}
This function is non-negative and non-decreasing. Moreover, one has $\Phi (x) = 0$ for 
$0 \leq x \leq h / ( 1 - |a|)$. Hence, since $| u (a) | \leq \frac{h}{1 - |a|}$\,, one has $\Phi \big( | u (a) | \big) = 0$.
\par
\smallskip

Then
\begin{align}
\int_h^2 \Phi '' \bigg( \frac{2h}{t}\bigg)\, \frac{h^2}{t^3}\, \tilde\nu (t)\,dt 
& = \int_{h}^{2 ( 1 - |a|)} \Phi '' \bigg( \frac{2h}{t}\bigg)\, \frac{h^2}{t^3}\, \tilde\nu (t)\,dt 
\label{morceau 2 bis} \\
& \leq \int_{h}^\infty \frac{h^2}{t^3}\,dt = \frac{1}{2}\,\cdot \notag
\end{align}
\smallskip

Now,
\begin{align*}
\Phi \Big( \frac{1}{2}\Big) 
& = \int_0^{1/2} \Phi ' (t) \,dt \geq \int_{1/4}^{1/2} \Phi ' (t) \,dt 
\geq \int_{1/4}^{1/2} \frac{t}{2} \, \Phi '' \Big( \frac{t}{2} \Big) \,dt \\
& \geq \Phi '' \Big( \frac{1}{8}\Big) \int_{1/4}^{1/2} \frac{t}{2} \,dt 
= \frac{3}{64} \, \Phi '' \Big( \frac{1}{8}\Big) \,.
\end{align*}
When $ h < (1 - |a|)/ 8$, one has $1/8 > h / (1 - |a|)$; hence $\Phi '' (1/8) = 1/ \tilde \nu (16 h)$, and 
$\Phi '' (1) = 1/ \tilde \nu (2h)$. We get hence, from \eqref{majoration avec Stanton bis}, \eqref{morceau 1}, 
\eqref{morceau 2} and \eqref{morceau 2 bis}:
\begin{equation}
m_\phi \big( S (1, h) \big) 
\leq \frac{64}{3}\, \tilde\nu (16 h) \bigg[ \frac{\tilde\nu (h)}{\tilde\nu (2h)} + 2 \bigg] 
\leq 64\, \tilde\nu (16h) \,.
\end{equation}
\par

Since $W (1, t) \subseteq S (1, 2 t)$, we get 
$m_\phi \big( W (1, h) \big) \leq 64\, \sup_{w \in S (1, 32 h)} N_\phi (w)$ for $0 < h < (1 - |\phi (0)|)/16$, 
and that ends the proof of Theorem~\ref{majoration mesure Carleson}, since $S (1, 32 h) \subseteq W (1, 64h)$.
\qed
\bigskip

\noindent{\bf Remark.} A slight modification of the proof gives the following improvement, if one allows 
a (much) bigger constant.

\begin{theoreme}\label{majo avec moyenne}
There are universal constants $C, c > 1$ such that 
\begin{displaymath}
m_\phi \big( S (\xi, h) \big) \leq 
C\, \frac{1}{A \big( S(\xi, c h) \big)} \, \int_{S(\xi, c h)} N_\phi (z)\,dA (z)
\end{displaymath}
for every analytic self-map $\phi \colon \D \to \D$, every $\xi \in \partial \D$, and 
$0 < h < (1 - |\phi (0)| )/ 8$.
\end{theoreme}

\noindent{\bf Proof.} We are going to follow the proof of Theorem~\ref{majoration mesure Carleson}. We 
shall assume that $\xi =1$ and we set:
\begin{equation}
I (t) = \int_{S (1, t)} N_\phi (z)\,dA (z) \,.
\end{equation}

Then:\par
1) When $|1  - z| < h$, , we have, instead of \eqref{morceau 1}:
\begin{align}
\int_{|1 - z| < h} \Phi '' \big( |u (z)| \big) \,|u '(z)|^2 N_\phi (z)  \,dA(z) 
& \leq \int_{S (1, h)} \Phi '' (1) \frac{1}{h^2}\,N_\phi (z)\, dA (z) \\
& = \Phi'' (1) \frac{1}{h^2}\, I(h)\,. \notag
\end{align}
\par

2) For $| z - 1 | \geq h$, we write:
\begin{align*}
\int_{|1 - z| \geq h} \Phi '' \big( |u (z)| \big) \,|u '(z)|^2 & N_\phi (z) \,dA(z) \\
& = \sum_{k=1}^\infty 
\int_{{}_{\scriptstyle kh \leq |1 - z| < (k + 1) h}} 
\hskip - 50 pt\Phi '' \big( |u (z)| \big) \,|u '(z)|^2 N_\phi (z) \,dA(z) \\
& \leq 4 \sum_{k=1}^\infty \Phi '' \Big( \frac{2h}{kh} \Big)\, \frac{h^2}{k^4 h^4}\, I\big( (k + 1) h \big) \\
& = 4 \sum_{k=1}^\infty \Phi '' \Big( \frac{2}{k} \Big) \, \frac{1}{k^4 h ^2} \, I \big( (k + 1) h \big) \,.
\end{align*}
We take, with $a = \phi (0)$:
\begin{equation}
\Phi '' (v) = \left\{
\begin{array}{cl}
\displaystyle 0 & \text{if } 0 \leq v \leq h / ( 1 - |a|) \,,\\
\smallskip
\displaystyle \frac{1}{I \big( (\frac{2}{v} + 1) h) \big)} & \text{if } v >  h / ( 1 - |a|) \,.
\end{array}
\right.
\end{equation}
Then
\begin{equation}
\int_{|1 - z| \geq h} \Phi '' \big( |u (z)| \big) \,|u '(z)|^2 N_\phi (z) \,dA(z) 
\leq \frac{4}{h^2} \sum_{k=1}^\infty \frac{1}{k^4} = \frac{4}{h^2} \, \frac{\pi^4}{90} \leq \frac{5}{h^2} 
\,\cdot
\end{equation}

Since $h < (1 - |a|)/8$, one has $1/8 > h/ (1 - |a|)$; hence $\Phi '' (1/8) = \frac{1}{I (17 h)}$ and 
$\Phi '' (1) = \frac{1}{I (3h)}$\,. Therefore:
\begin{align*}
m_\phi \big( S (1, h) \big) 
& \leq \frac{64}{3}\, I (17 h) \, \bigg[ \frac{1}{h^2}\, \frac{I (h)}{I (3h)} + \frac{5}{h^2} \bigg] \\
& \leq \frac{64}{3}\, I (17 h) \, \frac{6}{h^2} = 128 \frac{I (17 h)}{h^2} \\
& \leq 128 \times 17^2 \, \frac{I (17 h)}{A \big( S (1, 17 h) \big)}\, \raise 1pt \hbox{,}
\end{align*}
ending the proof of Theorem~\ref{majo avec moyenne}. \qed

 
\section{Some consequences}

In  \cite{CompOrli} (see also \cite{LLQR-CRAS}, Th\'eor\`eme~4.2), we proved (Theorem~4.19) that the Carleson 
function of an analytic self-map $\phi$ has the following property of homogeneity, improving that $m_\phi$ is a 
Carleson measure: $m_\phi \big( S (\xi, \eps\, h) \big) \leq K\,\eps\, m_\phi \big( S (\xi, h) \big)$ 
for $0 < h < 1 -|\phi (0)|$, $0 < \eps < 1$ and $\xi \in \partial \D$, where $K$ is a universal constant. It 
follows from Theorem~\ref{theo principal}, (actually Theorem~\ref{th majoration Nevanlinna} and 
Theorem~\ref{majoration mesure Carleson}) that:
\begin{theoreme}
There exist a universal constant $K > 0$ such that, for every analytic self-map $\phi$ of $\D$, 
one has, for $0 < \varepsilon < 1$:
\begin{equation}\label{homogeneite Nevanlinna}
\nu _\phi (\varepsilon \, t) \leq K\, \varepsilon\, \nu_\phi (t)\,, 
\end{equation}
for $t$ small enough.\par
More precisely, for $t$ small enough, one has, for every $\xi \in \partial \D$:
\begin{equation}
\nu_\phi (\xi, \eps\, t) \leq K\, \eps\, \nu_\phi (\xi, t) \,,
\end{equation}
where $\nu_\phi (\xi, s) = \sup_{w \in W (\xi, s) \cap \D} N_\phi (w)$.
\end{theoreme}

Note that the two above quoted theorems give Theorem~\ref{homogeneite Nevanlinna} {\it a priori} only for 
$0 < \eps < 1/K$; but if $1/K \leq \eps < 1$, one has 
$\nu_\phi (\xi, \eps\, t) \leq \nu_\phi (\xi, t) \leq K\, \eps\, \nu_\phi (\xi, t)$.
\bigskip

We shall end this paper with a consequences of Theorem~\ref{theo principal} for composition operators. 
Recall that if $\Psi$ is an Orlicz function, the Hardy-Orlicz space is the space of functions $f  \in H^1$ whose 
boundary values are in the Orlicz space $L^\Psi (\partial \D, m)$. 
We proved in \cite{CompOrli}, Theorem~4.18 (see also \cite{LLQR-CRAS}, Th\'eor\`eme~4.2) that, if  
$\frac{\Psi (x)}{x} \mathop{\longrightarrow}\limits_{x \to \infty} \infty$, the composition operator 
$C_\phi \colon H^\Psi \to H^\Psi$ is compact if and only if, for every $A > 0$, one has 
$\rho_\phi (h) = o\,\big[ 1/ \Psi \big( A \Psi^{-1} (1/h) \big) \big]$ when $h$ goes to $0$; in other words, if and 
only if 
\begin{displaymath}
\lim_{h \to 0} \frac {\Psi^{-1} (1/h)} {\Psi^{-1} \big(1/\rho_\phi (h) \big)} = 0\,.
\end{displaymath}
This remains true when $H^\Psi = H^1$. Hence Theorem~\ref{theo principal} gives:
\begin{theoreme}\label{Nevanlinna compact}
Let $\phi \colon \D \to \D$ be an analytic self-map and $\Psi$ be an Orlicz function. Then the 
composition operator $C_\phi \colon H^\Psi \to H^\Psi$ is compact if and only~if 
\begin{equation}\label{condition Nevanlinna compact}
\sup_{|w| \geq 1 - h} N_\phi (w) = o\,\bigg( \frac{1}{\Psi \big( A \Psi^{-1} (1/h) \big) } \bigg)\,, 
\quad \text{as } h \to 0\,, \qquad \forall A > 0.
\end{equation}
\end{theoreme}

It should be noted, due to the arbitrary $A > 0$, that \eqref{condition Nevanlinna compact} may be replaced by
\begin{equation}\label{condition bis Nevanlinna compact}
\sup_{|w| \geq 1 - h} N_\phi (w)  \leq\frac{1}{\Psi \big( A \Psi^{-1} (1/h) \big) } \,, 
\qquad \forall A > 0,
\end{equation}
for $h \leq h_A$, and this condition also writes, setting $\nu_\phi (h) =\sup_{|w| \geq 1 - h} N_\phi (w)$ 
(see \eqref{definition nu phi}):
\begin{equation}
\lim_{h \to 0} \frac{\Psi^{-1} (1/h)}{\Psi^{-1} \big( 1/ \nu_\phi (h) \big)} = 0\,.
\end{equation} 
\medskip

It is known that if $C_\phi\colon H^2 \to H^2$ is compact, then 
$\lim_{|z| \to 1} \frac{1 - |\phi (z)|}{1 - |z|} = \infty$, and that this condition is sufficient when 
$\phi$ is univalent, or finitely-valent, but not sufficient in general (see \cite{McCluer-Shapiro} and 
\cite{Shap-livre}, \S~3.2). It follows from Theorem~\ref{Nevanlinna compact} that an analogous 
result holds for Hardy-Orlicz spaces:
\begin{theoreme}
Let $\phi \colon \D \to \D$ be an analytic self-map, and $\Psi$ be an Orlicz function. Assume that the 
composition operator $C_\phi \colon H^\Psi \to H^\Psi$ is compact. Then:
\begin{equation}\label{condition compacite Nevanlinna}
\lim_{|z| \to 1} \frac{\Psi^{-1} \bigg(\displaystyle \frac{1}{1 - |z|} \bigg)} 
{\Psi^{-1} \bigg(\displaystyle  \frac{1}{1 - |\phi (z) |} \bigg)}  = \infty\,.
\end{equation}
\nobreak
Conversely, if $\phi$ is finitely-valent, then \eqref{condition compacite Nevanlinna} suffices for  
$C_\phi \colon H^\Psi \to H^\Psi$ to be compact.
\end{theoreme}

Recall that the assumption ``$\phi$ is finitely-valent'' means that there is an integer $p \geq 1$ such 
that each $w \in \phi (\D)$ is the image by $\phi$ of at most $p$ elements of $\D$.\par
\medskip

\noindent{\bf Proof.} To get the necessity, we could use Theorem~\ref{Nevanlinna compact} and the 
fact that $1 - |z| \leq \log \frac{1}{|z|} \leq N_\phi\big( \phi (z) \big)$; but we shall give a more 
elementary proof. 
Let $HM^\Psi$ be the closure of $H^\infty$ in $H^\Psi$. Since $C_\phi (H^\infty) \subseteq H^\infty$, 
$C_\phi$ maps $HM^\Psi$ into itself and $C_\phi \colon H^\Psi \to H^\Psi$ being compact, its restriction 
$C_\phi \colon HM^\Psi \to HM^\Psi$ is compact too. We know that the evaluation 
$\delta_a \colon f\in HM^\Psi \mapsto f (a) \in \C$ has norm $\approx \Psi^{-1} \big( \frac{1}{1 - |a|} \big)$ 
(\cite{CompOrli}, Lemma~3.11); hence $\delta_a/ \| \delta_a\| \mathop{\longrightarrow}\limits_{|a| \to 1} 0$ 
weak-star (because $|\delta_a (f) | = |f (a) | \leq \|f \|_\infty$ for $f \in H^\infty$). If $C_\phi$ is compact, 
its adjoint $C_\phi^\ast$ also; we get hence 
$\| C_\phi^\ast (\delta_a / \| \delta_a\|) \| \mathop{\longrightarrow}\limits_{|a| \to 1} 0$. 
But $C_\phi^\ast \delta_a = \delta_{\phi (a)}$. Therefore 
\begin{displaymath}
\frac{\Psi^{-1} \bigg( \displaystyle \frac{1}{1 - |\phi (a) |} \bigg)} 
{\Psi^{-1} \bigg( \displaystyle \frac{1}{1 - |a|} \bigg) } \mathop{\longrightarrow}_{|a| \to 1} 0\,.
\end{displaymath}
\par

Conversely, assume that \eqref{condition compacite Nevanlinna} holds. For every $A > 0$, one 
has, for $|z|$ close enough to $1$: 
$\Psi^{-1} \big(\frac{1}{1 - |z|} \big) \geq A\,\Psi^{-1} \big( \frac{1}{1 - |\phi (z) |} \big)$; in other words, 
one has: 
$1 / \Psi \big( A \Psi^{-1} (1/ 1 - |\phi (z)|) \big) \geq 1 - |z|$. But, when $\phi$ is $p$-valent, 
and if $w = \phi (z)$ with $|z| > 0$ minimal, one has $N_\phi (w) \leq p \log \frac{1}{|z|} \approx 1 -|z|$. 
Since $|z| \to 1$ when $|w| = |\phi (z)| \to 1$ (otherwise, we should have a sequence $(z_n)$ converging 
to some $z_0 \in \D$ and $\phi (z_n)$ would converge to $\phi (z_0) \in \D$), we get 
$\sup_{|w| \geq 1 - h} N_\phi (w) \lesssim 1 / \Psi \big( A \Psi^{-1} (1/ 1 - |w|) \big) 
\leq 1 / \Psi \big( A \Psi^{-1} (1/ 1 - h) \big) $, for $h$ small enough. By Theorem~\ref{Nevanlinna compact},  
with \eqref{condition bis Nevanlinna compact}, that means that $C_\phi$ is compact on $H^\Psi$.
\hfill $\square$
\bigskip

Other consequences will be given in the subsequent paper~\cite{LLQR_09}.
\bigskip

\noindent{\bf Acknowledgement.} This work has been initiated when the second-named author was invitated 
by the Departamento de An\'alisis Matem\'atico of the Universidad de Sevilla, in April 2007; it is a pleasure 
to thanks all its members for their warm hospitality. Part of this work was also made during the fourth-named 
author visited the University of Lille 1 and the University of Artois (Lens) in June 2007. The fourth-named 
author is partially supported by a Spanish research project MTM2006-05622.


\bigskip

\vbox{\noindent{\it 
{\rm Pascal Lef\`evre}, Univ Lille Nord de France F-59\kern 1mm 000 LILLE, \\
U-Artois, Laboratoire de Math\'ematiques de Lens EA~2462, \\
F\'ed\'eration CNRS Nord-Pas-de-Calais FR~2956, \\
F-62\kern 1mm 300 LENS,
FRANCE \\ 
pascal.lefevre@euler.univ-artois.fr 
\smallskip

\noindent
{\rm Daniel Li}, Univ Lille Nord de France F-59\kern 1mm 000 LILLE, \\
U-Artois, Laboratoire de Math\'ematiques de Lens EA~2462, \\
F\'ed\'eration CNRS Nord-Pas-de-Calais FR~2956, \\
F-62\kern 1mm 300 LENS,
Facult\'e des Sciences Jean Perrin,\\
Rue Jean Souvraz, S.P.\kern 1mm 18, FRANCE \\ 
daniel.li@euler.univ-artois.fr
\smallskip

\noindent
{\rm Herv\'e Queff\'elec}, Univ Lille Nord de France F-59\kern 1mm 000 LILLE, \\
USTL, Laboratoire Paul Painlev\'e U.M.R. CNRS 8524, \\
F-59\kern 1mm 655 VILLENEUVE D'ASCQ Cedex, 
FRANCE \\ 
queff@math.univ-lille1.fr
\smallskip

\noindent
{\rm Luis Rodr{\'\i}guez-Piazza}, Universidad de Sevilla, \\
Facultad de Matem\'aticas, Departamento de An\'alisis Matem\'atico,\\ 
Apartado de Correos 1160,\\
41\kern 1mm 080 SEVILLA, SPAIN \\ 
piazza@us.es\par}
}

\end{document}